\documentclass[12pt]{article}

\usepackage{times}
\usepackage{a4}
\usepackage{amsmath,amsfonts,bm}
\usepackage{harvard}
\usepackage{graphicx}

\newtheorem{theorem}{Theorem}
\newtheorem{lemma}[theorem]{Lemma}

\newcommand{\BI}{%
\begin{list}{$\bullet$}{\itemsep=0pt\parsep=4pt\listparindent=1em\topsep=-4pt}}
\newcommand{\EI}{\end{list}}

\begin{document}

\begin{titlepage}

\baselineskip=18pt

\begin{center}

{\bf \Large Models and algorithms for skip-free Markov decision processes on trees}
 
\vspace{2cm}
 
{\large E.J. Collins\\
Department of Mathematics,\\ University of Bristol,\\
 University Walk,\\ Bristol BS8 1TW, UK.}
 
\vspace{1cm}

Accepted for publication subject to minor changes by the Journal of the Operational Research Society (JORS)

\vspace{1cm}

\begin{quote}
We introduce a class of models for multidimensional control problems
which we call {\em skip-free Markov decision processes on trees}.
We describe and analyse an algorithm applicable to Markov decision 
processes of this type that are skip-free in the negative direction.
Starting with the finite average cost case, 
we show that the algorithm combines the advantages of 
both value iteration and policy iteration 
-- it is guaranteed to converge to an optimal policy and 
optimal value function after a finite number of iterations 
but the computational effort required for each iteration step 
is comparable with that for value iteration.
We show that the algorithm can also be used to solve 
discounted cost models and continuous time models, and
that a suitably modified algorithm can be used to solve communicating models.

\end{quote}

\end{center}

\vfill

\noindent {\small KEYWORDS: 
Multidimensional Markov decision processes, Dynamic programming, Queueing, Inventory, Maintenance, Reliability
\\}

\end{titlepage}

\section{Introduction} \label{sec:intro}

Markov decision processes (MDPs) provide a class of stochastic optimisation models
that have found wide applicability to problems in Operational Research.
The standard methods for computing an optimal policy 
are based on value iteration, policy iteration and linear programming algorithms
\cite{Whi93}.
Each approach has its advantages and disadvantages.
In particular, each step in value iteration is relatively computationally 
inexpensive but the value function may take some time to converge 
and the algorithm provides no direct check that it has computed
the optimal value function and an optimal policy.
Conversely, each step in policy iteration may be computationally expensive
but the algorithm can be proved to converge in a finite number of steps,
confirms when it has converged and automatically identifies the optimal 
value function and an optimal policy on exit.

Here we focus on models with special structure, in that they
are {\em skip-free in the negative direction} \cite[p.10]{Kei65}
or {\em skip-free to the left} \cite{StWe89}; 
i.e.\ whatever the action taken,
the process cannot pass from one state to a `lower'
state without passing through all the intervening states. 
Such skip-free models arise naturally in many areas where OR is applied. 
The most obvious examples are the control of discrete time random walks 
and continuous time birth and death processes \cite{Ser81}
such as queueing control problems with single unit arrivals and departures
(see, for example, \citeasnoun{StWe89} and references therein). 
In these basic one-dimensional models, the state space $S$ is (a subset of) the integer 
lattice and transitions are only possible to the next higher or lower integer state.
However there are several other standard OR models that fall within the wider
one-dimensional skip-free framework, including examples from the areas of
queueing control with batch arrivals \cite{StWe89},
inventory control \cite{Mil81} and reliability and maintenance \cite{Der70,Tho82}.

Previous treatments of controlled skip-free processes
have considered only the one-dimensional formulation.
For processes with the `skip-free to the left' property, work has focused 
on qualitative properties, in particular the existence of monotone optimal 
policies for models with appropriately structured cost functions \cite{StWe89,StWe99}.
Conversely, work on processes with the corresponding `skip-free to the right' 
property has concentrated on analysis of an approximating bisection method 
for countable state space models \cite{WiSt86,WiSt00}.
We note that skip-free type ideas have also been exploited in a different direction by 
\cite{Whi05} and citing authors, 
where the emphasis has been on reducing the computational complexity
associated with policy iteration for quasi birth-death processes.

An intuitive way of characterising the essential features of our finite 
skip-free recurrent model is that 
the model is skip-free if and only if the state space 
can be identified with the graph of a finite tree, rooted at $0$,
with each state $i$ corresponding to a unique node in the tree,
and such that for every action $a \in A$,
the only possible transitions from state $i$ under action $a$ are either to 
its `parent' state  or to a state in the subtree rooted at $i$,
with appropriate modifications for state $0$ which has no parent
and for terminal nodes which have only a parent and no descendants.

In this setting, the one-dimensional skip-free model above, with state space 
$S = \{0,1,\ldots,M\}$, corresponds to the simplest case where 
the tree reduces to a single linearly ordered branch connecting the 
root node $0$ through states $1, 2, \ldots, M-1$ to the terminal node $M$,
and transitions from state $i$ are possible only to states $ j \in \{i-1, i, \ldots,M \}$. 
However, the analysis extends easily to cases with a 
richer, possibly multidimensional, state space, 
where the appropriate model is in terms of transitions on a finite tree. 
Examples of genuinely skip-free models with multidimensional state spaces
arise in simple multi-class queueing systems with batch arrivals
\cite[and references therein]{YeSe94,He00},
but such treatments have focused mainly on describing the behaviour
of the process for a fixed set of parameters (actions) rather than 
comparing actions in an optimality framework. 

The rest of the paper is organized as follows. 
We start by describing models for average cost finite state recurrent MDPs 
that are skip-free in the negative direction,
illustrating our approach with a motivating example.
We then propose a skip-free algorithm 
that combines the advantages of values iteration and policy iteration:
the computational effort required for each iteration step 
is comparable with that for value iteration,
but the algorithm is guaranteed to converge after a finite number of iterations
and automatically identifies the optimal value function and an optimal policy on exit.
We go on to show that the algorithm can be also be used to solve 
discounted cost models and continuous time models, and
that a suitably modified algorithm can be used to solve communicating models.
Finally, we build on the relationship between 
the average cost problem and a corresponding 
$x$-{\em revised first passage problem} 
to provide a proof of the main theorem and
identify other possible variants of the algorithm.


\section{The skip-free MDP model} \label{sec:sf-model}

Consider a discrete time Markov decision process (MDP) 
with finite state space $S$
over an infinite time horizon $t \in \{0,1,2,\ldots \}$.
Associated with each state $i \in S$ is a non-empty finite set of 
possible actions; 
since $S$ is finite, we assume without loss of generality
that the set of actions $A$ is the same for each $i$.
If action $a \in A$ is chosen when the process is in state $X_t = i$
at time $t$,
then the process incurs an immediate cost $c_i(a)$
and the next state is $X_{t+1} = j$ with probability $p_{ij}(a)$.

A policy $\pi$ is a sequence of (possibly history dependent and randomised) 
rules for choosing the action at each given time point $t$.
A {\em deterministic} decision rule corresponds to a function $d\!:\!S\!\rightarrow\!A$ 
and specifies taking action $a = d(i)$ when the process is in state $i$.
A {\em stationary deterministic} policy is one which 
always uses same the deterministic decision rule at each time point $t$.
Where the meaning is clear from the context, we use the same notation 
$d$ for both the decision rule and the corresponding stationary deterministic policy.

The expected average cost incurred by a policy $\pi$ with initial state $i$ is given by
$g_{\pi}(i) = \limsup_{n \to \infty} 
\frac{1}{n} \; E_{\pi} \left( \sum_{t=0}^{n-1} c_{X_t}(a_t) | X_0 = i \right ),$
where $X_t$ is the state at time $t$ and $a_t$ is the action chosen at time $t$ under $\pi$.
Similarly, for a given discount factor $0 < \beta < 1$,
the total expected discounted cost incurred by a policy $\pi$ with initial state $i$ 
is given by
$V^{\beta}_{\pi}(i) = E_{\pi} \left( \sum_{t=0}^{\infty} 
\beta^n \, c_{X_t}(a_t) | X_0 = i \right ).$


We say an MDP model is {\em recurrent} if the transition matrix corresponding to
every stationary deterministic policy consists of a single recurrent class.
We  say an MDP model is {\em communicating} if, for every pair of states $i$ and $j$ in $S$, 
$j$ is reachable from $i$ under some (stationary deterministic) policy $d$;
i.e.\ there exists a policy $d$, with corresponding transition matrix $P_d$,
and an integer $\, n \ge 0$, such that $P_d(X_n = j | X_0 = i) > 0$. 

When $S = \{ 0, 1, 2, \ldots, M \}$ is a subset of the integer lattice,
we say the MDP model is {\em skip-free in the negative direction}
\cite{Kei65,StWe89} if 
$p_{ij}(a) = 0$ for all $j < i-1$ and $a \in A$,
i.e.\ the process cannot move from state $i$ to a state with index $j < i$ 
without passing through all the intermediate states.
We will often find it easier to work in terms of the upper tail probabilities
$\bar{p}_{i j}(a) \equiv P(X_{t+1} \ge j \, | \, X_t = i, A_t = a) = \sum_{s=j}^M p_{i s}(a)$.
To avoid degeneracy, we assume that $p_{00}(a) < 1$ for $a \in A$
and that for each $i \in \{1, \ldots,M \}$, $p_{i i-1}(a) > 0$ for at least one $a \in A$.
In this setting, a recurrent model requires that, for all $a \in A$, 
$p_{ii-1}(a) > 0$ for $i = 1,\ldots,M$ and $p_{ii}(a) < 1$ for all $i \in S$.
In contrast a communicating model allows 
there to be $i$ and $a$ with $p_{ii-1}(a) = 0$ and /or $p_{ii}(a) = 1$.

To apply this idea in a wider context, we note that the essence of a skip free model is that:
~(i) there is a single distinguished state, say $0$; 
~(ii) for any other state $i$ there is a unique shortest path from $i$ to $0$;
~(iii) from each state $i \ne 0$ the process can only make transitions to either 
the adjacent state in the unique path from $i$ to $0$,
or to some state $j$ for which $i$ lies in the unique shortest path from $j$ to $0$.

In the finite one dimensional case, for each $k$ 
there is exactly one state for which the shortest path to state $0$ has length $k$.
Thus there is a $1$--$1$ mapping of the states to the integers
$\{0, 1, \ldots,M\}$ such that the distinguished state maps to $0$
and the state for which the shortest path had length $k$ maps to $k$.
In a more general setting, for each $k$  
there may be more than one state for which the shortest path has length $k$.
In this case, rather than $S$ mapping to the integer lattice,
there is a fixed tree ${\cal T}$ (in the graph theoretic sense)
such that each state corresponds to a unique node of the tree, 
with the distinguished state mapping to the root node.
It may help to visualise movement between states in terms
of the corresponding movement between nodes on the tree.


To formalise this general model, we start by considering a finite rooted tree ${\cal T}$ 
with $N + 1$ nodes labelled $0,1,2,\dots,N$, with root node $0$, and with a given edge set.
The tree structure implies that for each pair of nodes $i$ and $j$ 
there is a unique minimal path (set of edges) in the tree that connects $i$ and $j$.
Thus the nodes in the tree can be partitioned into level sets $L_0 = \{0\}, L_1, \ldots, L_M$ 
such that, for $m = 0, \ldots,M-1$, $i \in L_{m+1}$ if and only if 
the minimal path from $i$ to $0$ passes through exactly $m$ intermediate nodes.

For adjacent nodes $i \in L_m$ and $j \in L_{m+1}$,
we say $i$ is the parent of $j$ and $j$ is a child of $i$
if the minimal path from $j$ to $0$ passes through $i$. 
More generally, for $i \in L_m$ and $j \in L_r, \, r > m$,
we say $j$ is a descendant of $i$ if the minimal path from $j$ to $0$
passes through $i$. Each node $j \ne 0$ has a unique parent. 
We write $\rho(j)$ for the parent of $j$,
we write ${\cal D}(j)$ for the set of descendants of $j$,
and we write ${\cal T}(j) \subset {\cal T}$ for (the nodes of the) 
sub-tree rooted at $j$, so ${\cal T}(j) = \{j\} \cup {\cal D}(j)$.
A state with no descendants is said to be a terminal state,
so all states in the highest level $L_M$ are terminal states.
For simplicity of presentation we will assume that these are the only terminal states;
the analysis easily extends to cases where intermediate levels $L_m$ can also contain some terminal states.
For each $j \in {\cal D}(i)$, we write $\Delta(i,j)$ for the set of states 
following $i$ in the unique minimal path in the tree connecting  $i$ to $j$,
so if the path passes through $s-1$ intermediate states 
and takes the form $i = r_0 \rightarrow r_1 \rightarrow \cdots \rightarrow r_s = j$,
then  $\Delta(i,j) = \{r_1,\ldots,r_s\}$.

Now consider a finite MDP with state space $S$ and action space $A$.
Assume we can construct a rooted tree ${\cal T}$ such that 
(i) the states in $S$ correspond to the nodes of ${\cal T}$, 
and (ii) for every state $i \in S$ and action $a \in A$,
the only possible transitions from state $i$ under action $a$ are either to 
its parent state $\rho(i)$
or to a state in the subtree ${\cal T}(i)$ rooted at $i$,
with appropriate modifications for state $0$ which has no parent
and for terminal nodes which have only a parent and no descendants.
We will say that such an MDP is {\em skip-free (in the negative direction) on the tree } ${\cal T}$.
As with the integer lattice model above,
it is often convenient work in terms of the the upper tail probabilities
$\bar{p}_{ij}(a) = P(X_{t+1}\in {\cal T}(j)|X_t=i, A_t=a),$
corresponding to the probability that the next transition from state $i$
under action $a$ is to a state in the subtree rooted at $j$.


To illustrate and motivate the general case, where a multidimensional model is required, 
consider (\cite{He00,YeSe94}) a single-server multi-class queueing system  
with $K > 1 $ customer classes and finite capacity $M$ (including the job, 
if any, in service). 
Assume the service discipline is pre-emptive but otherwise takes no account of class. 
A job that arrives when the system is not full enters service immediately 
and the job currently in service at that point returns to the head of the buffer. 
When a job completes service, the server next serves the job at the head of the buffer. 
Any job that arrives when the system is full is lost.

The model is most naturally formulated in continuous time, 
with exponential inter-arrival and service time distributions, 
though it can easily be translated to a discrete time setting 
using the methods of section \ref{sec:cont}. 
Assume class $k$ jobs arrive at rate $\lambda_k$ 
and complete service at class and action dependent rate $\mu_k(a)$, 
where different actions $a \in A$ correspond to different service levels. 
Since the model needs to keep track of the class of each job as it enters service, 
we take the state to be the multidimensional vector 
$\bm{i} = (i_1, \ldots, i_{M})$ 
where $i_1$ denotes the class of the job currently in service, 
$i_m$ denotes the class of the job waiting for service 
in the buffer in place $m, \; m = 2, \ldots, M$, 
and $i_m = 0$ if the $m$th place is empty. 
Assume costs are incurred at rate $c(\bm{i},a)$ reflecting both holding costs and action costs.

The possible transitions under the model are 
the completion of the job currently in service,
corresponding to the transition
$\bm{i} = (i_1, \ldots, i_{M}) \rightarrow
(i_2, \ldots, i_{M},0)$,
or the arrival of a class $k$ job ($k = 1,\ldots,K$) 
to a partially full system,
corresponding to the transition 
$\bm{i} = (i_1, \ldots, i_{M}) \rightarrow
\bm{j} = (k,i_1, \ldots,i_{M-1})$. 

For $M \ge 2$ this model cannot be represented as a skip-free MDP 
with linear structure, i.e.\ with each state $\bm{i}$
having exactly one child $\bm{j}$ with $\bm{i} = \rho(\bm{j})$.
To see this, 
let $\bm{a}$ denote the state $(a, i_2, \ldots, i_{M})$ with $i_M \ne 0$, 
let $\bm{b}$ denote the state $(b, i_2, \ldots, i_{M})$,
differing from $\bm{a}$ in only the first component, 
and let $\bm{c}$ denote the state $(i_2, \ldots, i_{M},0)$.
The only possible direct transitions to and from $\bm{a}$ are from $\bm{c}$ and to $\bm{c}$.
Similarly for $\bm{b}$. 
If $\bm{c}$ is restricted to having just one child, then the only possibilities
are either  
(i) $\bm{a}$ has no parent (so $\bm{a}$ is the root state), 
$a = \rho(\bm{c})$ and $c = \rho(\bm{b})$, or
(ii) $\bm{b}$ has no parent (so $\bm{b}$ is the root state), 
$b = \rho(\bm{c})$ and $c = \rho(\bm{a})$. 
In case (i), $\bm{b}$ can have no children 
so none of the other states can reach the root state
as they cannot reach $\bm{b}$ in a skip-free manner under any policy;
in case (ii) $\bm{a}$ can have no children and a similar argument applies.


\vspace*{-0.2in}
\begin{figure}[ht] 
\resizebox{12cm}{8cm}{\includegraphics{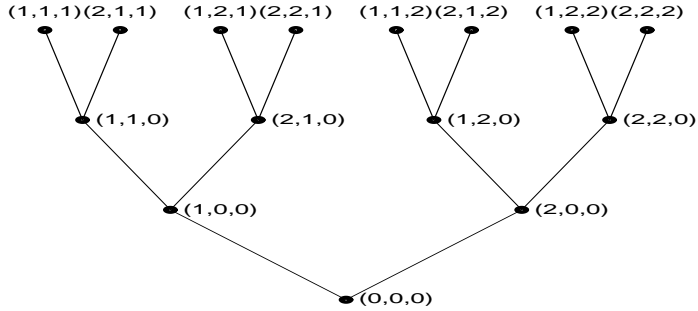}}
\vspace*{-0.7in}
\caption{The tree ${\cal T}$ corresponding to 
the state space for the pre-emptive multi-class queueing system 
of with $K=2$ job classes and capacity $M = 3$.} \label{fig:pe}
\end{figure}


However we can represent the model as a skip-free MDP on a tree ${\cal T}$ as follows. 
We take $L_0 = \{ (0, \ldots,0) \}$ to contain the state corresponding to the empty queue 
and take the level sets $L_m, \; m=1, \ldots,M$ to each contain the $K^m$ states 
of the form $\bm{i} = (i_1, \ldots,i_m,0,\ldots,0)$. 
Given a state $\bm{i} = (i_1, \ldots,i_M) \in L_m$  
we assign it parent $\rho(\bm{i}) = (i_2, \ldots, i_M, 0)$ 
and assign it $K$ children of the form $\bm{j} = (k,i_1, \ldots,i_{M-1}), \; k = 1, \ldots,K$ 
(with appropriate modifications for $L_0$ and $L_M$).
The set of descendants ${\cal D}(\bm{i})$ is the set of all states 
of the form $(k_1, \ldots,k_r, i_1, \ldots,i_m, 0, \ldots,0)$ for $r = 1, \ldots, M-m$ 
(where there are $M-m-r$ trailing $0$s).
The possible transitions under the model
correspond exactly to transitions
from $\bm{i}$ to its parent $\rho(\bm{i})$
or to one of its $K$ children, so 
the MDP satsifies the conditions required for it to be skip free
in the negative direction on the tree ${\cal T}$.
Figure \ref{fig:pe} illustrates the tree corresponding to the 
state space for a system with $K=2$ job classes and capacity $M = 3$.
Extensions with direct transitions to more general descendants,
of form $(k, \ldots,k, i_1, \ldots,i_m, 0, \ldots,0)$
are possible if batch arrivals are allowed, subject to appropriate capacity constraints.
  

\section{The skip-free algorithm} \label{sec:sf-alg}

For finite recurrent MDP models, the solution to the 
expected average cost problem can be characterised by the corresponding 
{\em average cost optimality equations} \cite[\S 8.4]{Put94} 
\begin{align}
h_i &= \min\nolimits_{a \in A} \{ \; c_i(a) - g + \sum\nolimits_{j \in S} p_{ij}(a) h_j \; \}   
& i \in S \label{eq:oe}
\end{align} 
\noindent in that 
(i) there exist real numbers $g^*$ and $h^*_i, i \in S$ satisfying the optimality 
equations;
(ii) the optimal average cost is the same for each initial state and is given by $g^*$;
(iii) the optimality equations uniquely determine $g^*$ and determine the $h^*_i$ 
up to an arbitrary additive constant;
(iv)
the stationary deterministic policy $d^*$ is an average cost optimal policy, 
where, for each $i \in S$, $d^*(i)$ is an action achieving $
\min_a \{ \; c_i(a)+\sum_{j \in S} p_{ij}(a)h^*_j \, \}.$

It follows from (iv) above that there is an optimal policy in the class of 
stationary deterministic policies. 
We therefore restrict attention from now on to stationary deterministic policies,
writing `policy' as a shorthand for `stationary deterministic policy'
and writing $g(d)$ for the average cost under a given stationary deterministic policy $d$.

For each $i,j \in S$, we can interpret $h^*_i - h^*_j$ as 
the asymptotic relative difference in the total cost
that results from starting the process in state $i$ rather than state $j$, 
under the stationary deterministic policy $d^*$. 
Thus the quantities $h^*_i - h^*_j$ are uniquely defined,
but the quantities $h^*_i, i \in S$ are defined only up to an arbitrary additive constant.
We focus on the particular solution normalised by setting $h^*_0=0$
and refer to the corresponding $h^*_i$ as the normalised relative costs under an optimal 
policy.


In general, the optimality equations~(\ref{eq:oe}) cannot be solved directly. 
Instead an optimal policy in the class of stationary deterministic policies is
usually found by methods based on value iteration, policy iteration or linear programming,
or combinations of these approaches \cite{Put94}.
For skip-free models, however, we have the following simplification. 

\begin{lemma} \label{lem:oe}
For finite recurrent skip-free average cost MDPs, 
the optimality equations (\ref{eq:oe}) are equivalent to the equations 
\begin{subequations} \label{eq:oe-sfall}
\begin{align} 
y_i &= \min\nolimits_a \{ \; (c_i(a) - x)/p_{i \rho(i)}(a) \; \} \; & i \in L_M  \label{eq:oe-sfa} \\
y_i &= \min\nolimits_a \{ \; (c_i(a) - x + \sum\nolimits_{k \in {\cal D}(i)} \bar{p}_{i k}(a) y_k )/p_{i \rho(i)}(a) 
\; \}  & i \in L_{M-1},\ldots,L_1 \label{eq:oe-sfb} \\
0  &= \min\nolimits_a \{ \; c_0(a) - x + \sum\nolimits_{k \in {\cal D}(0)}  \bar{p}_{0 k}(a) y_k \; \} \label{eq:oe-sfc}
\end{align}
\end{subequations}
\noindent in that 
(i) these equations also have unique solutions $x$ and $y_i, \, i \in {\cal D}(0)$;
(ii) the optimal average cost is $g^* = x$ and the normalised relative costs under an optimal
policy satisfy $h^*_i - h^*_{\rho(i)} = y_i, \; i \in {\cal D}(0)$;
(iii) an optimal stationary deterministic policy is given by $d^*$, 
where $d^*(i)$ is any action minimising the rhs of the corresponding equation for $y_i$
and $a_0$ is an action minimising the rhs in (\ref{eq:oe-sfc}).
\end{lemma}
\noindent {\bf Proof}
For skip-free models, 
the only possible transitions from state $i \in {\cal D}(0)$
are to state $\rho(i)$, to state $i$ itself, or to a state $j \in {\cal D}(i)$.
Thus equations~(\ref{eq:oe}) take the form
\begin{align} 
h_i  &= \min\nolimits_{a \in A} \{ \; c_i(a) - g + \sum\nolimits_{j \in {\cal D}(i)} p_{i j}(a) h_j + p_{i i}(a) h_{i} 
+ p_{i \rho(i)}(a) h_{\rho(i)} \; \}  & i \in S \label{eq:new-oe}
\end{align}
with appropriate modification to give the normalised solution with $h_0 = 0$.
Values $h_i$ and $g$ satisfy  (\ref{eq:new-oe}) 
if and only if in each equation $h_i \le $ 
the rhs for all $a$, with equality for at least one $a$.
With appropriate modifications for the root node $0$ and for terminal nodes,
simple rearrangement in shows that 
$h_i \le c_i(a) - g + \sum_{j \in {\cal D}(i)} p_{i j}(a) h_j + p_{i i}(a) h_{i} 
+ p_{i \rho(i)}(a) h_{\rho(i)}$
if and only if 
$p_{i \rho(i)}(a) (h_i - h_{\rho(i)}) \le  c_i(a) - g + \sum_{j \in {\cal D}(i)} p_{i j}(a) (h_j - h_i),$
and that equality in one expression implies equality in the other. 

Now write $x$ for $g$ and for each $i \in {\cal D}(0)$ write $y_i$ for $h_i - h_{\rho(i)}$.
For each $j \ne i \in {\cal D}(i)$, write $\Delta(i,j) = \{r_1,\ldots,r_s\}$ 
for the states following $i$ in the unique minimal path from $j$ to $i$.
For each $k = 1,\ldots,s$, $r_{k-1}$ is the parent of $r_k$ so that $r_{k-1} = \rho(r_k)$.
Hence 
$h_j - h_i
= h_{r_s} - h_{r_0}
= \sum_{k=1}^s (h_{r_k} - h_{r_{k-1}})
= \sum_{k=1}^s (h_{r_k} - h_{\rho(r_k)})
= \sum_{k=1}^s y_{r_k}
= \sum_{r \in \Delta(i,j)} y_{r}$.
Now if $j $ is a descendant of $i$ 
and $r \ne j$ is in the path connecting  $i$ and $j$,
then $r$ is a descendant of $i$ and $j$ is in the subtree rooted at $r$,
and vice versa. 
Thus for fixed $i$ and $a$ we have that 
$\sum_{j \in {\cal D}(i)} p_{i j}(a) (h_j - h_i)
= \sum_{j \in {\cal D}(i)} \sum_{r \in \Delta(i,j)} p_{i j}(a) y_{r} 
= \sum_{r \in {\cal D}(i)} \sum_{j \in {\cal T}(r)} p_{i j}(a) y_r 
= \sum_{r \in {\cal D}(i)}  \bar{p}_{ir}(a) y_r$. 

Taking account of the modifications for the root state $i=0$ 
and the terminal states $i \in L_M$,
and the fact that $i \in L_m \implies {\cal D}(i) \subset L_{m+1} \cup \cdots \cup L_M$,
it follows that there are $g$ and $h_i$ satisfying (\ref{eq:new-oe}) if and only if
there are values $x$ and $y_i$ satisfying (\ref{eq:oe-sfall}).
  \hfill $\Box$
  

In the optimality equations ~(\ref{eq:oe-sfall}), 
the value of $y_i, \, i \in L_M$ depends only on $x$,
and in each subsequent equation the value of $y_i$ depends only
on $x$ and the values of $y_k$ for $k \in {\cal D}(i)$.
Thus, if the value of $x$ was known, it would be easy
to compute the $y_i$ in turn for $y_i \in L_M, \ldots,L_1$
and to determine the corresponding policy which takes 
the optimal action in each state $i \in S$.

This observation motivates an iterative approach to finding an average cost optimal policy:
(i) choose an initial policy $d_0$ and compute its expected average cost $g_0 = g(d_0)$;
(ii) given a current policy $d_n$ with expected average cost $g_n$,
compute an updated policy $d_{n+1}$ by setting $x = g_n$
and solving (\ref{eq:oe-sfa}) and (\ref{eq:oe-sfb}), 
and compute its expected average cost $g_{n+1}$;
(iii) iterate until convergence.
This approach forms the basis for the following {\em skip-free} algorithm.
Its properties are set out in the subsequent theorem.\\

\noindent {\bf Skip-free algorithm}

\noindent 1. \underline{Initialisation}: 

\noindent Choose an arbitrary  initial policy $d_0$.
Perform a single iteration of step 2 below, with $x = 0$ 
and with $a_i$ restricted to the single value $d_0(i)$, $i \in S$.
Set $g_0 = u_0$.

\medskip

\noindent 2. \underline{Iteration}: 

\noindent Set $x  = g_n$. 

\medskip

\noindent $\bullet$ For $i \in L_M$ compute:\\  
$\begin{array}{llcl}
& a_i & = & \mathrm{argmin}_a \{ \; (c_i(a) - x)/p_{i \rho(i)}(a) \; \}  \\
& y_i & = & (c_i(a_i) - x)/p_{i \rho(i)}(a_i)  \\
& t_i & = & 1/p_{i \rho(i)}(a_i) \\
\end{array}$

\medskip

\noindent $\bullet$ For $i \in L_r, \, r =  M-1,\ldots,1$ compute:\\ 
$\begin{array}{llcl}
& a_i & = & \mathrm{argmin}_a \{ \; (c_i(a) - x + \sum_{k \in {\cal D}(i)} \bar{p}_{i k}(a) y_k )
/p_{i \rho(i)}(a) \; \} \\
& y_i & = & (c_i(a_i) - x + \sum_{k \in {\cal D}(i)} \bar{p}_{i k}(a_i) y_k )/p_{i \rho(i)}(a_i)\\
& t_i & = &  (1 + \sum_{k \in {\cal D}(i)} \bar{p}_{i k}(a_i) )/p_{i \rho(i)}(a_i)\\
\end{array}$

\medskip

\noindent $\bullet$ For $j = 0$ compute:\\
$\begin{array}{llcl}
& a_0 & = & \mathrm{argmin}_a \{ \; (c_0(a) - x + 
\sum_{k \in {\cal D}(0)}   \bar{p}_{0 k}(a) y_k) /
            (1+ \sum_{k \in {\cal D}(0)}   \bar{p}_{0 k}(a_0) t_k) \; \} \\
& u_0 & = & (c_0(a_0) - x + \sum_{k \in {\cal D}(0)}   \bar{p}_{0 k}(a_0) y_k) / (1+ \sum_{k \in {\cal D}(0)}   \bar{p}_{0 k}(a_0) t_k) \\
& t_0 & = & (1+ \sum_{k \in {\cal D}(0)}   \bar{p}_{0 k}(a_0) t_k) / (1 - p_{0 0}(a_0)) \\

\end{array}$

\medskip

\noindent Set $d_{n+1}(i) = a_i$ for $i \in S$ 
and set $g_{n+1} = g_n + u_0 $.

\medskip

\noindent 3. \underline{Termination}:

\noindent If $u_0 < 0$ then return to step 2. 

\medskip

\noindent If $u_0 = 0$ then stop. 
Return $d_{n+1}$ as an optimal policy, 
return $g_{n+1}$ as the optimal average cost,
and for each $i \in {\cal D}(0)$ return $h_i = \sum_{j \in \Delta(0,i)}y_j$
as the corresponding normalised relative cost.

\begin{theorem} \label{th:grr}
Consider the skip-free algorithm above applied to a finite recurrent
skip-free average cost MDP model. Then:\\
(i) At each iteration either
$g_{n+1} < g_n$, so $d_{n+1}$ is a strict improvement on $d_n$,
or $g_{n+1} = g_n$.
In the latter case $g_{n+1} = g^*$, $d_{n+1}$ is an optimal average cost policy,
and the corresponding normalised relative costs 
are given by $h_0^* =0, \, h_j^* = \sum_{i \in \Delta(0,j)}y_i, \; j \in {\cal D}(0)$.\\
(ii) The algorithm converges after a finite number of iterations.
\end{theorem}

\noindent {\bf Remarks}
(1) The motivation for the particular choice of action in state $0$ is given in 
the remarks following the proof of the theorem.
(2) The updates are particularly simple in the one dimensional case 
where $S = \{ 0,1,\ldots,M\}$. 
Here $\sum_{k \in {\cal D}(i)}$ simplifies to $\sum_{k =i+1}^M$
and $\rho(i)$ simplifies to  $i-1$.
(3) The computational requirement for each iteration in step 2 of the algorithm
is clearly similar to that of the corresponding step in value iteration, in that it only
requires simple  evaluations rather than the solution of a set of equations. 
While the algorithm is also similar to policy evaluation in that it returns the average cost 
of policy $d_n$ at the end on the $n$th iteration, 
it differs from standard policy iteration in that
it the values of $y_i$ returned do not correspond to the relative costs 
under $d_n$. Only at convergence do the relative costs and average cost correspond 
to the same (optimal) policy.
(4) The basic principle underlying this iterative approach appears to be
similar to that used in \cite{Low74}, but 
the results there were restricted to a very specific model 
with simple birth and death structure.
Other treatments of skip-free models \cite{WiSt86,StWe89,StWe99,WiSt00}
have used iterative methods to search for a good approximation for the average cost $x$,
based on the value of current and previous approximations,
or used the form of the optimality equations to 
derive qualitative properties of the solution,
in particular monotonicity of optimal policies,
but neither approach explicitly identified the 
simple skip-free improvement algorithm described here. 


\section{Discounted, continuous and communicating models} \label{sec:var}

The skip-free algorithm can also be used to solve discounted cost  
and continuous time problems, in each case by transforming the problem
into an equivalent average cost problem. 
Moreover, a suitably modified algorithm can be used to solve communicating models.
For ease of presentation, we focus on the one dimensional case,
indicating how the argument can be extended to the general model as required.


\subsection{Discounted cost models} \label{sec:disc}

Consider a recurrent  MDP model that is skip-free in the negative direction,
with state space $S = \{0, 1, \ldots,M\}$,
finite action space $A$, 
transition probabilities  $p_{ij}(a)$, immediate costs $c_i(a)$
and discount factor $\beta$.
Following \citeasnoun[p.31]{Der70}, we construct an average cost MDP 
with modified state space $\{0,1, \ldots,M,M+1 \}$
and modified transition probabilities and immediate costs given by:
\begin{align}
&p'_{ij}(a)        =  \beta p_{ij}(a),   &&c'_i(a)     =  c_i(a), & i,j = 0,1, \ldots,M, 
\; \; a \in A \nonumber\\
&p'_{M+1 \, M}(a)  = \beta,  \quad       &&c_{M+1}(a)  =  0,      & a \in A \nonumber\\
&p'_{i \, M+1}(a)  = 1 - \beta,          &&                       & i = 0,1, \ldots,M + 1, 
\; \; a \in A \nonumber
\end{align} 

\noindent In the spirit of similar models \cite{Low74,WiSt86}, we note that
this new average cost MDP inherits from the original model the property of being 
skip-free in the negative direction.

Let $g'$ and $h'_i, \, i=0,\ldots,M+1$ be the optimal average cost and the corresponding 
relative costs 
for the new average cost problem, normalised by setting $h'_0=0$.
From above, $g'$ and $h'_i, \, i=1,\ldots,M+1$, are the unique solutions 
to the optimality equations~(\ref{eq:oe}), 
and any set of actions achieving the minimum on the rhs defines an optimal policy.
In terms of the original parameters, these equations take the form
\begin{align}
&h'_{M+1} &&=  - g' + \beta h'_M + (1-\beta) h'_{M+1} \nonumber\\
&h'_i     &&=   \min_a \{ \, c_i(a) - g' + \beta \sum_{j=0}^{M} p_{ij}(a) h'_j + 
(1 - \beta)h'_{M+1} \, \} & i = 0,\ldots,M \nonumber
\end{align}
Now set $v_j = h'_j - h'_{M+1} + g'/(1-\beta), \; j = 0,\ldots,M$.
Then rewriting the equations for $h_0, \ldots, h_M$ in terms of $v_0,\ldots,v_M$,
we see that the $v_i$ satisfy the equations
\begin{align}
v_i  &=  \min_a \{ \, c_i(a)  + \beta \sum_{j=0}^{M} p_{ij}(a) v_j \, \} &i = 0,\ldots,M. \nonumber
\end{align} 
Thus the $v_j$ satisfy the optimality equations for the discounted cost problem,
and so represent the unique optimal $\beta$ discounted cost function \cite[p.148]{Put94}.

Finally, let $x'$ and $y'_0, \ldots,y'_{M+1}$ be solutions to the policy iteration algorithm 
applied to the new skip-free average cost problem.
Then $g' = x'$ and  $h'_j =  y'_j + \cdots + y'_1, \, j=1, \ldots,M + 1$. 
Thus the optimal value function for the discounted problem is given 
explicitly in terms of the output of the policy iteration algorithm by
\begin{align}
v_j &= x'/(1-\beta) -(y'_{j+1} + \cdots + y'_{M+1}) & j = 0,\ldots,M \nonumber
\end{align}
and a policy which is optimal for the modified average cost problem 
is also optimal for the original discounted cost problem.

The extension to the general skip-free MDP tree model is straightforward,
requiring just the addition of an extra state for each terminal state (node) 
to preserve the skip-free property. This extra state now becomes 
the terminal node in that branch. 
Transitions from this extra state are to the corresponding previous terminal
node, with probability $\beta$, or back to itself, with probability $1-\beta$.
Transition probabilities from non-terminal states are modified as above, 
by setting $p'_{ij}(a) = \beta p_{ij}(a)$ if $j$ is a non-terminal 
node of the modified tree and by assigning the remaining transition probability 
$1-\beta$ to the newly added terminal nodes of the modified sub-tree ${\cal T}(i)$ rooted at $i$. 
The precise assignment may be chosen arbitrarily -- for example, each new 
terminal node in the modified sub-tree may be chosen with equal probability
-- as long as the total probability sums to $1 - \beta$.


\subsection{Continuous time models} \label{sec:cont}

Consider a continuous time Markov decision process (CTMDP)
with finite state space $S$ and finite action space $A$.
Assume that when the current action is $a$ and the process is in 
state $X_t = i$, the process incurs costs at rate $c_i(a)$
and makes transitions to state $j \in S$ at rate $q_{ij}(a)$
(where transitions back to the same state are allowed).
For infinite horizon problems, under either an average cost or 
a discounted cost criterion,
we can restrict attention to stationary policies and to models in 
which decisions are made only at transition epochs \cite[p.560]{Put94}.
For simplicity of presentation we again restrict attention to recurrent 
models and defer treatment of unichain and communicating models to
Section \ref{sec:com}. 
As for MDPs, we say a CTMDP is skip-free in the negative direction 
if the process cannot move from each state $i$ to a state
$j < i$ without passing through all the intermediate states, 
i.e.\ $q_{ij}(a) = 0$ for all $j < i-1$ and $a \in A$.

To apply the skip-free algorithm, we first convert the model to an equivalent uniformised model
\cite{Lip75} with rate $\Lambda = \max_{i \in S \; a \in A} \sum_{j \in S} q_{ij}(a)$.
In this model, when the current action is $a$ and the process is in state $i$,
transitions back to state $i$ occur at rate
$\Lambda - \sum_{j \ne i} q_{ij}(a)$ while transitions 
to state $j \ne i$ occur at rate $q_{ij}(a)$,
so that overall transitions occur at uniform rate $\Lambda$.
Next we construct a discrete time problem with the same state and action space, 
where for $i,j \in S$ and $a \in A$ the transition probabilities and immediate costs are given by
$p'_{ij}(a) =  q_{ij}(a)/\Lambda,\; i \ne j; p'_{ii}(a) = 1 - \sum_{j \ne i} q_{ij}(a)/\Lambda; 
c'_i(a) =  \Lambda c_i(a)$. 
If the original CTMDP is recurrent and skip-free,
then the discretised model is recurrent and skip-free 
and can be solved using the algorithm.

Finally, let $d'$ and $g'$ be the optimal policy and the optimal average cost
identified by the algorithm for the discrete time problem. 
Then the optimal policy $d^*$ and the optimal average cost $g^*$
for the uniformised continuous time problem are the same 
as $d'$ and $g'$, 
and the normalised relative costs for the uniformised problem
are given in terms of those for the discrete problem
by $h_i^* = h'_i/\Lambda , \; i \in S$ \cite[\S 11.5]{Put94}. 


\subsection{Communicating models} \label{sec:com}

So far we have assumed the MDP model is recurrent. 
There are natural applications for which this assumption excludes sensible policies,
such as policies that are recurrent only on a strict subset of $S$.
Simple examples include: 
maintenance/replacement problems where a policy might specify replacing an item when 
the state reached some lower level $K > 0$ with a item of level $L < M$;
inventory problems where a policy might reorder when the stock reached some 
lower level $K > 0$ and/or reorder up to level $L < M$;
queueing control problems where a policy might turn the server off when the queue size 
reached some lower level $K > 0$ and/or might refuse to admit new entrants when the queue 
size reached level $L < M$.
In each case, determining optimal values for $K$ and $L$ might be part of the problem.
In this section we extend our result to the wider class of communicating MDP models,
to enable us to address examples like these. 

We  say an MDP model is {\em communicating} if, for every pair of states $i$ and $j$ in $S$, 
$j$ is reachable from $i$ under some (stationary deterministic) policy $d$;
i.e.\ there exists a policy $d$, with corresponding transition matrix $P_d$,
and an integer $\, n \ge 0$, such that $P_d(X_n = j | X_0 = i) > 0$.
We say that $d$ is {\em unichain} if it decomposes $S$ into a single recurrent class
plus a (possibly empty) set of transient states;
if there is more than one recurrent class we say $d$ is {\em multichain}.
Let $d$ be a multichain policy and, for each $k$,
let $g_k$ denote the average cost under $d$ starting in a state in $E_k$,
and let $E_m$ be a recurrent set with smallest average cost, say $g_m$. 
Because the model is skip-free, $E_m$ must consist of a sequence of 
consecutive states $K_m,\ldots,L_m$; 
again, because the model is skip-free, the action in each each state $j$ greater than $L_m$
can be changed if necessary so that $E_m$ is reachable from $j$;
finally, because the model is communicating, the action in each state $j$ less than $K_m$
can be changed if necessary so that $E_m$ is reachable from $j$.
Denote by $d'$ the new policy created by changing actions in this way, 
if necessary, but leaving the actions in $E_m$ unchanged. 
Then $d'$ is unichain by construction, 
and the average cost starting in each state $j \in S$ 
is $g_m$, which is no greater than the average cost starting in $j$ under $d$.
Thus, for average cost skip-free communicating models, 
nothing is lost by restricting attention to unichain policies.

In contrast to recurrent models, communicating models allow 
there to be $i$ and $a$ with $p_{ii}(a) = 1$ and/or $p_{ii-1}(a) = 0$.
For each $r = 0,1,\ldots,M$, 
let $U_r$ be the (possibly empty) set of unichain policies $d$ for which 
$p_{rr-1}(d(r)) = 0$ but $p_{ii-1}(d(i)) > 0$ for $i = r+1, \ldots,M$
(where we take $p_{ii-1}(a) \equiv 0$ for all $a$ for $i=0$).
Every unichain policy must be in $U_r$ for some $r$.
Partition the possible actions for each state $i \in S$ into 
$B_i = \{a \in A : p_{ii-1}(a) >0 \}$
and its complement $\bar{B}_i = \{a \in A : p_{ii-1}(a) = 0 \}$,
where $\bar{B}_i$ may be empty but 
$B_i$ is non-empty by the assumptions of the skip free model in Section~\ref{sec:sf-model}.
Then for a unichain policy $d \in U_r$, we have that
$d(i) \in B_i, \; i = r+1, \ldots,M$; 
that state $r$ is recurrent and $d(r) \in \bar{B}_r$ by definition;
and that states $i < r$ are transient.

Thus the minimum average cost over policies in $U_r$
is the same as the minimum average cost for a modified skip-free MDP model $\Pi_r$
with the same transition probabilities and immediate costs but 
with reduced state space $S_r = \{r,\ldots,M\}$
and with state-dependent action spaces 
$A_i = B_i$ for $i = r+1, \ldots,M$ and $A_r = \bar{B}_r$.
In this notation, the model of Section~\ref{sec:sf-model} corresponds to $\Pi_0$
and state $r$ plays the same role as the recurrent distinguished 
state in $\Pi_r$ that state $0$ plays in $\Pi_0$.
If we compare the result of applying the skip-free algorithm to $\Pi_r$ 
with the result of applying it to $\Pi_0$,
we see that, for the same current value of $x$,
the algorithm computes the same values of 
$y_i$, $t_i$, and $a_i$ in states $i = M,M-1,\ldots,r+1$.
However, in state $r$, the skip-free algorithm applied to $\Pi_r$ 
computes quantities appropriate to the distinguished state,
say $a^r$ and $u^r$, where\\

$\begin{array}{lcl}
a^r & = & \mathrm{argmin}_{a \in \bar{B}_r}
            \{ \; (c_r(a) - x + \sum_{k=r+1}^M   \bar{p}_{r k}(a) y_k) /
            (1+ \sum_{k=r+1}^M   \bar{p}_{r k}(a) t_k) \; \} \\
u^r & = & (c_r(a^r) - x + \sum_{k=r+1}^M   \bar{p}_{r k}(a^r) y_k) / 
            (1+ \sum_{k=r+1}^M   \bar{p}_{r k}(a^r) t_k) 
\end{array}$\\

\noindent
and computes an updated `minimising' policy $d^r_{n+1}$ 
with average cost $g^r_{n+1}$,where \\

$\begin{array}{llcl}
 d^r_{n+1}(r) & = & a^r; \; \; d^r_{n+1}(i) = a_i, \; i=r+1,\ldots,M, \quad \mbox{and}\\
 g^r_{n+1} & = & x + u^r.
\end{array}$\\

This motivates the following modified skip-free algorithm. 
First, it includes these extra  computations for each state $r$,
so that, in a single iteration, it simultaneously computes
the optimal policy $d^r_{n+1}$ and its average cost $g^r_{n+1}$ 
for each $S_r$.
Secondly, at the end of the $n-1$th iteration 
it sets $x = g_{n} = \min_r g^r_{n}$, and sets $d_{n}$ to be the corresponding policy,
where ties are broken by choosing the $d^r_{n}$ with the smallest index $r$.
Say the minimum average cost at this stage is achieved by a policy with index $r = K$
Then, by the properties of the skip-free algorithm applied to $\Pi_K$,
at the end of the next iteration 
either (i) $g^K_{n+1} < g^K_{n} = x$,
in which case $g_{n+1} = \min_r g^r_{n+1} < x = g_n$; 
or (ii) $u^K_{n+1}=0$ and $g^K_{n+1} = g^K_{n} = x = \min_r g^r_{n+1}$,
so $g_{n+1} = g_n$ and $d_{n+1} = d^K_{n+1}$ is an optimal average cost policy
for starting states $i = K, \ldots,M$.
In this case, because the model is communicating,
it is possible \cite[p.351]{Put94} to modify the actions chosen by the policy
in the, now transient, states $0, \ldots, K-1$ so that the modified $d_{n+1}$ satisfies the 
optimality equations for all states $0, \ldots, M$
and is an average cost optimal policy.
We summarise this discussion in the following theorem.

\begin{theorem} 
Consider the skip-free algorithm modified as above applied to a finite communicating discrete time
average cost skip-free MDP model with state space $S = \{0,1,2,\ldots,M\}$. Then:\\
(i) At each iteration of the skip-free algorithm either
$g_{n+1} < g_n$ and $d_{n+1}$ is a strict improvement on $d_n$,
or $g_{n+1} = g_n$ and for some $K$ 
the policy satisfies the optimality equations for states $K, \ldots, M$.\\
(ii) The modified skip-free algorithm converges after a finite number of iterations.
\end{theorem}


Finally, note that it is easy to check if a skip-free model is communicating.
An assumption of the (non-degenerate) skip-free model was that
each state $i<M$ was reachable from $i+1$.
It follows that a skip-free MDP with state space $S = \{0,1,\ldots,M\}$ 
is communicating if and only if $M$ is reachable from $0$ under at least one stationary 
deterministic policy $d$.
Let $N_0 = 0$,  let $N_1$ be the index of the maximum state $j$ for which $p_{0j}(a) > 0$ 
for some $a \in A$, and for $m = 1, 2, \dots $ let $N_{m+1}$ be the index of the maximum 
state $j$ for which $p_{ij}(a) > 0$ for some $0 \le  i \le N_m$ and $a \in A$. 
As the state space is finite, the sequence $\{N_m\}$ terminates, say with state $N$.
Since the model is skip-free, $N$ is the largest state that is reachable by all states below 
it, and the model is communicating if and only if $N = M$.

The extension to a general skip-free communicating models is straightforward. 
Again, the idea is that for each state $i$ the skip-free algorithm
is modified so that in passing it solves the corresponding sub-problem 
$\Pi_i$ with state space ${\cal T}(i)$ and with state $i$ as the distinguished state,
and then computes the optimal updated average cost and policy
by minimising over the costs and policies for each of the sub-problems. 


\section{Proof of Theorem \ref{th:grr}} \label{sec:proofs}

We start our analysis of the average cost MDP model by defining a 
related problem (or class of problems) 
that we will call the $x$-{\em revised first return problem}.
The model for this problem has the same state space $S$, 
the same action space $A$ and the same 
transition probabilities $\{p_{ij}(a)\}$ as the average cost model.
However, for each fixed $x$, the immediate costs in the corresponding 
$x$-revised problem are revised downward by $x$, 
so $c_i(a)$ is revised to $c_i(a) - x$.
Whereas the original problem was to find a policy $d$ 
that minimised the expected average cost $g(d)$, 
the objective for this new problem is to find a policy that minimises
the expected $x$-revised cost until first return to state $0$,
where, for a process starting with $X_0 = 0$, we define the first return epoch 
to state $0$ to be the smallest value $\tau > 0$ such that $X_{\tau-1} \ne 0$ 
and $X_\tau = 0$. The MDP is assumed recurrent under any stationary
deterministic policy, so $\tau$ is well defined and almost surely finite.

For a fixed policy $d$, starting in state $0$, write 
$\tau(d)$ for the expected first return epoch under $d$, 
$C(d)$ for the expected first return cost under $d$, and 
$H(d,x)$ for the expected $x$-revised first return cost under $d$.
The average costs and the $x$-revised costs under $d$
are related by the equations
\begin{equation} \label{eq:hct}
g(d) = C(d)/\tau(d),  \quad H(d,x) = C(d) - x \tau(d), \quad g(d) = x + H(d,x)/\tau(d),
\end{equation} 
where the first equation follows from viewing the average cost problem 
from a renewal-reward perspective \cite[p.160]{Ros70}
and noting that state $0$ is recurrent under any stationary deterministic policy $d$,
and the second follows from noting that the expected $x$-revised cost under $d$ until first return to state $0$
is just the original expected cost $C(d)$ adjusted downwards by an amount $x$ for 
an expected time period $\tau(d)$.  

\begin{lemma} \label{lem:dfr}
For fixed $x$, let $a_i, \, i \in {\cal D}(0)$ be actions minimising 
the rhs in equations~(\ref{eq:oe-sfa}) and (\ref{eq:oe-sfb})
and let $y_i, \, i \in {\cal D}(0)$ be the corresponding $y$ values.
Set
\begin{equation} \label{eq:a-fr}
a_{0} = \mathrm{argmin}_a \{ \; (c_0(a) - x + \sum\nolimits_{k \in {\cal D}(0)}
\bar{p}_{0 k}(a) y_k) / (1 - p_{0 0}(a)) \; \}.
\end{equation}
and let $d$ be the policy that takes action $a_i$ in state $i, \; i \in S$.
Then $d$  minimises the expected $x$-revised cost until first return to state $0$,
and the expected $x$-revised first return cost under  $d$ is
\begin{equation} \label{eq:h-fr}
H(d, x) = (c_0(a_{0}) - x + \sum\nolimits_{k \in {\cal D}(0)}  \bar{p}_{0 k}(a_{0}) y_k) / (1 - p_{0 0}(a_{0})).
\end{equation}
\end{lemma}

\noindent {\bf Proof}
Since the process is Markov and skip-free in the negative direction, 
it follows that a policy minimises the expected $x$-revised cost 
until first return to state $0$ if and only if it also minimises the expected $x$-revised 
total cost until first passage to state $0$ for each starting state 
$i \ne 0 \in $, i.e.\ $i \in {\cal D}(0)$,
 and hence minimises the expected cost until first passage from $i$ to to $\rho(i)$
for each $i \in {\cal D}(0) $. 
For the one-dimensional case where $S = \{0, 1, \ldots,M\}$, 
this problem has been called the $x$-revised first passage problem \cite{StWe89}.
For fixed $x$ and $i \in \{1, \ldots,M \}$, let $a_i$ be actions minimising 
the rhs in equations~(\ref{eq:oe-sfa}) and (\ref{eq:oe-sfb}) 
and let $y_i$ be the corresponding $y$ values.
Then they show that the policy $d$ that takes action $d(i) = a_i$ in state $i$
is optimal for the $x$-revised first passage problem
and the minimal expected cost until first passage from $i$ to $i-1$ is given by $y_i$.
With only minor notational changes, their results extend directly to the general 
case where $S$ corresponds to the nodes of a tree, 
$\{ 1, \ldots,M\}$ is replaced by ${\cal D}(0)$
and $i-1$ is replaced by $\rho(i)$.
It follows that the policy that uses actions $a_i$ in $i \in {\cal D}(0) $
has the property that for each state $i$ it also minimises 
the expected total $x$-revised cost until first passage to state $0$ and that 
the minimum expected $x$-revised total cost until first passage 
to state $0$, starting in state $i \ne 0$, is given by the sum of the $y_i$
values along the path from $i$ to $0$, 
i.e.\ $\sum\nolimits_{k \in \Delta(0,i)} y_k$. 

Now consider a process that starts in state $0$.
Under a policy that specifies action $a$ in state $0$,
the expected time until the process first leaves state $0$ is $1/(1-p_{0 0}(a))$ and
during that time it incurs $x$-revised costs at rate $c_0(a) - x$ per unit time.
Conditional on leaving state $0$, the first transition is to state $j$ 
with probability $p_{0 j}(a)/(1-p_{0 0}(a))$.
From above, the minimum additional expected total cost until the process next
re-enters state $0$ is $\sum\nolimits_{k \in \Delta(0,j)} y_k$,
and this minimum expected cost is achieved by the policy that takes actions 
$a_i$ in states $i \in {\cal D}(0) $.
Thus, if a policy $d$ takes action $a$ in state $0$,
the minimum expected $x$-revised cost from leaving state $0$ 
until first return to state $0$ is 
$ H(d,x) = \sum\nolimits_{j \in {\cal D}(0)}  p_{0 j}(a) \sum\nolimits_{k \in \Delta(0,j)} y_k /(1 - p_{0 0}(a))$
$=$$ \sum\nolimits_{k \in {\cal D}(0)} \sum\nolimits_{j \in {\cal T}(k)} p_{0 j}(a) y_k /(1 - p_{0 0}(a))$
$=$$ \sum\nolimits_{k \in {\cal D}(0)} \bar{p}_{0 k}(a) y_k /(1 - p_{0 0}(a))$.
It follows that the optimal action in state $0$ is one that minimises the quantity
$(c_0(a) - x + \sum\nolimits_{k \in {\cal D}(0)} \bar{p}_{0 k}(a) y_k) / (1 - p_{0 0}(a))$
and the expected $x$-revised first return cost $H(d,x)$ is as shown. \hfill $\Box$

\begin{lemma} \label{lem:imp}
Let $d$ be a fixed policy with expected average cost $g(d)$
and let $d^1$ be the optimal $x$-revised policy specified in Lemma \ref{lem:dfr} 
for the case $x = g(d)$.
Then:\\
(i) the average cost under $d^1$ is no greater than the average cost under $d$,\\
(ii) if the average cost under $d^1$ is the same as the average cost under $d$
then $d^1$ is an optimal policy for the average cost problem.
\end{lemma}

\noindent {\bf Proof}
(i) For the fixed $x$, we know from Lemma \ref{lem:dfr} that 
$d^1$ is an optimal policy  for the $x$-revised first return problem. 
Thus $H(d^1,x)) \le H(d,x)$, and from (\ref{eq:hct}) this implies 
$C(d^1) - x \tau(d^1) \le C(d) - x \tau(d)$.
Because $x$ corresponds to the average cost under $d$,
then, from (\ref{eq:hct}), $x = g(d) = C(d)/\tau(d)$ so $C(d) - x \tau(d) = 0$.
Thus, $H(d^1) = C(d^1) - x \tau(d^1) \le 0$ 
and $g(d^1) = C(d^1)/\tau(d^1) \le x = g(d)$.\\
(ii)
If $g(d^1) = g(d)$, then from above $H(d^1,x) = H(d) = 0$.
But, from Lemma~\ref{lem:dfr},
$H(d^1,x) = (c_0(a_0) - x + \sum_{k=1}^M  
\bar{p}_{0 k}(a_0) y_k) / (1 - p_{0 0}(a_0))$,
where $p_{0 0}(a_0) < 1$.
It follows that  $H(d^1,x) = 0 \implies 
(c_0(a_0) - x + \sum_{k=1}^M \bar{p}_{0 k}(a_0) y_k) = 0$. 
Thus, when $g(d^1) = g(d)$,
the values $x = g(d^1)$ and the corresponding values of $y_i, \, i \in {\cal D}(0)$ satisfy 
the optimality equations~(\ref{eq:oe-sfa}-\ref{eq:oe-sfc})
and $d^1$ is a decision rule corresponding to the actions minmising the 
rhs of each equation.
It follows that $d^1$ is an optimal average cost policy, 
the optimal average cost is $g^* = g(d^1) = g(d)$
and the normalised relative costs under the optimal policy are 
$h^*_j = \sum_{k \in \Delta(0,j)}y_k$. \hfill $\Box$

\begin{lemma} \label{lem:eval}
Let $a_i, \, i \in S$ be fixed actions
and let $d$ be the fixed policy for which $d(i) = a_i, \, i \in S$.
Perform a single iteration of step 2 of the skip-free algorithm 
with starting value $x$ and with the action in each state $i$ restricted to the single value $a_i$.
If the algorithm output values are $u_0$, $y_i, \, i \in {\cal D}(0)$ and $t_i, \, i \in S$, 
then $H(d,x)$ and $\tau(d)$ are given by equations (\ref{eq:h-fr})
and (\ref{eq:timec}).
Further, if the starting value is $x = 0$, then $g(d) = u_0$.
\end{lemma}
\noindent {\bf Proof} 
The expression for $H(d,x)$ follows from Lemma \ref{lem:dfr} 
by considering the possible actions in state $i$ to be 
restricted to just the given $a_i$. 

For the expected first return epoch under $d$, write $t_0 = \tau(d) > 0$ 
and write $t_i > 0$ for the expected first passage time $t_i$ from $i$ to $i-1$.
Interpret $t_0$ as the expected $0$-revised first return cost under $d$
for a model with immediate costs $c_i(a) = 1$ for all states and actions (and with $x = 0$), 
with a similar interpretation for the $t_i$.
Then, as with the $y_i$, the $t_i$ can be computed recursively 
using the equations 
$t_i = 1/p_{i \rho(i)}(a_i), \, i \in L_M; 
t_i = (1 + \sum_{k \in{\cal D}(i)}  \bar{p}_{i k}(a_i) t_k )/p_{i \rho(i)}(a_i), 
\; i \in  L_{M-1},\ldots,L_1,$ and 
\begin{equation} \label{eq:timec} 
\tau(d) = t_0 = (1 + \sum\nolimits_{k \in {\cal D}(0)}    \bar{p}_{0 k}(a_0) t_k) / (1 - p_{0 0}(a_0)). 
\end{equation}

Finally set $x = 0$. Then $g(d) = H(d,0)/\tau(d)$ from  (\ref{eq:hct}), 
so from (\ref{eq:h-fr}) and (\ref{eq:timec}) 
$ g(d) = (c_0(a_0) + \sum\nolimits_{k \in {\cal D}(0)} \bar{p}_{0 k}(a_0) y_k) /
(1 + \sum\nolimits_{k \in {\cal D}(0)}    \bar{p}_{0 k}(a_0) t_k) = u_0.$ \hfill $\Box$


Given a current policy $d$ with average cost $x = g(d)$, 
both the original optimality equations (\ref{eq:oe-sfall})
and the $x$-revised approach suggest updating $d$ with a 
policy that for $i \in {\cal D}(0)$ uses the actions $a_i$ identified by
equations (\ref{eq:oe-sfa}) and (\ref{eq:oe-sfb}).
However they differ in their suggested action $a_0$ in state $0$ --
the former suggests using the action minimising 
the rhs in equation~(\ref{eq:oe-sfc}) while the latter suggests using 
the action identified in (\ref{eq:a-fr}). However, the above lemma suggests 
another possible choice would be   
\begin{equation} \label{eq:a-mi}
a_{0} = \mathrm{argmin}_a \{ \; (c_0(a) - x + \sum\nolimits_{k \in {\cal D}(0)}    
\bar{p}_{0 k}(a) y_k) / (1+ \sum\nolimits_{k \in {\cal D}(0)}    \bar{p}_{0 k}(a_0) t_k) \; \}.
\end{equation}
This results in a policy that minimises the average cost over all policies that take 
the given actions $a_i$ in states $i \in {\cal D}(0)$.
The next lemma shows all three variations either strictly improve on $d$ 
or identify an optimal policy.


\begin{lemma} \label{lem:123}
Let $d$ be a fixed policy and let $x = g(d)$.
For this $x$, let $a_i, \, i \in {\cal D}(0)$ be actions minimising 
the rhs in equations~(\ref{eq:oe-sfa}) and (\ref{eq:oe-sfb})
and let $y_i, \, i \in {\cal D}(0)$ be the corresponding $y$ values.
Let $a_0^1$ be the action specified by equation (\ref{eq:a-fr}),  
let $a_0^2$ be the action minimising the rhs of equation (\ref{eq:oe-sfc}),
and let $a_0^3$ be the action specified by equation (\ref{eq:a-mi}). 
For $k = 1, 2, 3$, let $d^k$ be the policy that takes action $a_i$ in state $i \in {\cal D}(0)$
and takes action $a_0^k$ in state $0$.
Then either (i) all three policies $d^k$ satisfy $g(d^k) < g(d)$, or  
(ii) all three policies satisfy $g(d^k)  = g(d)$ and each of the three (and $d$ itself) provides an optimal average cost policy.
\end{lemma}

\noindent {\bf Proof}
For fixed $x$ and any policy $d$, $g(d) - x = H(d,x)/\tau(d)$ from  (\ref{eq:hct}) 
and $\tau(d)$ is positive, so $g(d) - x$ has the same sign as $H(d,x)$.
Since all three policies take actions $a_i$ in states $i \in {\cal D}(0)$,
expression (\ref{eq:h-fr}) gives their respective expected $x$-revised first return 
costs as
$H(d^k,x) = (c_0(a_0^k) - x + \sum\nolimits_{k \in {\cal D}(0)}  \bar{p}_{0 k}(a_0^k) y_k) / (1 - p_{0 0}(a_0^k))$,
where each $p_{00}(a_0^k) < 1$ by the assumptions of the skip-free model. 

Now $H(d^2,x) < 0 \implies
(c_0(a_0^2) - x + \sum\nolimits_{k \in {\cal D}(0)}     \bar{p}_{0 k}(a_0^2) y_k) / (1 - p_{0 0}(a_0^2)) < 0 \implies 
(c_0(a_0^1) - x + \sum\nolimits_{k \in {\cal D}(0)}    \bar{p}_{0 k}(a_0^1) y_k) / (1 - p_{0 0}(a_0^1)) < 0$  
(as $a_0^1$ minimises this quantity over choice of $a$)
$\implies  H(d^1,x)< 0$.
Conversely $H(d^1,x) < 0 \implies
(c_0(a_0^1) - x + \sum\nolimits_{k \in {\cal D}(0)}     \bar{p}_{0 k}(a_0^1) y_k) / (1 - p_{0 0}(a_0^1)) < 0 \implies 
(c_0(a_0^1) - x + \sum\nolimits_{k \in {\cal D}(0)}     \bar{p}_{0 k}(a_0^1) y_k) < 0 \implies  
(c_0(a_0^2) - x + \sum\nolimits_{k \in {\cal D}(0)}     \bar{p}_{0 k}(a_0^2) y_k) < 0$  
(as $a_0^2$ minimises this quantity over choice of $a$)
$\implies  H(d^2,x)< 0$.
A similar argument utilising the definition of $a_0^3$
and the positivity of $(1+ \sum\nolimits_{k \in {\cal D}(0)}     \bar{p}_{0 k}(a_0) t_k)$
shows that $H(d^2,x) < 0 \Longleftrightarrow H(d^3,x) < 0$.
Exactly similar arguments then show
that $H(d^1,x) = 0 \Longleftrightarrow H(d^2,x) = 0 \Longleftrightarrow H(d^3,x) = 0$,
and that $H(d^1,x) > 0 \Longleftrightarrow H(d^2,x) > 0 
\Longleftrightarrow H(d^3,x) > 0$. 
The second part of the lemma then follows from Lemma~\ref{lem:imp}. \hfill $\Box$


\noindent {\bf Proof of Theorem \ref{th:grr}}
(i) It follows from Lemma \ref{lem:eval} that the initialisation step outputs $g_{0} = g(d_{0})$.
Now let $x = g_n$ and assume $g_n = g(d_n)$.
Then iteration $n+1$ outputs $g_{n+1} = g_n + u_0$,
where $u_0 = 
(c_0(a_0) - x + \sum\nolimits_{k \in {\cal D}(0)}   \bar{p}_{0 k}(a_0) y_k) / (1+ \sum\nolimits_{k \in {\cal D}(0)}    \bar{p}_{0 k}(a_0) t_k) = 
H(d_{n+1},x)/\tau(d_{n+1})$ from (\ref{eq:h-fr}) and (\ref{eq:timec}).
Thus $g_{n+1} = x + H(d_{n+1},x)/\tau(d_{n+1}) = g(d_{n+1})$ from equation (\ref{eq:hct}).
Since $g_0 = g(d_0)$, it follows by induction that $g_n = g(d_n)$ for $n = 0,1,2,\ldots$.

By construction at iteration $n+1$ the skip-free algorithm 
specifies $d_{n+1}(i) = a_i, \; i \in S$, 
where $a_i, \, i \in {\cal D}(0)$ are the actions minimising 
the rhs in equations~(\ref{eq:oe-sfa}) and (\ref{eq:oe-sfb}) for this value of $x$
(and $y_i, \, i \in {\cal D}(0)$ and $t_i, \, i \in {\cal D}(0)$ are the corresponding $y$ and $t$ values),
and $a_0$ is the action minimising the rhs in equation (\ref{eq:a-mi}).
It follows from Lemma \ref{lem:123} that either 
$g(d_{n+1}) < g(d_n)$, or $g(d_{n+1}) = g(d_n)$ and both $d_{n+1}$ and $d_n$  provide optimal average cost policies.
Finally the expression for $h_j^*$ follows from considering the case $i = 0$
in the representation $h_j - h_i =  \sum_{k \in \Delta(i,j)}y_k$ in Lemma \ref{lem:oe}
with the normalisation $h_0 = 0$.
 
(ii) Since the set of possible stationary deterministic decision rules is finite,
and each iteration prior to convergence leads to a strict improvement
and hence a strictly different decision rule, 
the process must converge after a finite number of steps. \hfill $\Box$


\noindent{\bf Remark}
The update proposed in the skip-free algorithm uses $a_0$ 
satisfying (\ref{eq:a-mi}). 
It has the property that, for each current policy $d$,
it generates an improved policy with average cost at least as small as 
the other two variants considered in Lemma \ref{lem:123}.
This does not guarantee that improvements
using this update converge faster than improvements using either
of the other two variants. After one iteration, each policy
may generate a different starting point for the next iteration,
and our results do not allow us to compare the policies from these different 
starting points -- indeed it might be that the larger the improvement from the 
first iteration, the smaller the improvement resulting from the second iteration,
as the average cost is now closer to the optimal value.
Our experience has been that the number of iterations taken by all three methods was 
often the same. 
Where one was fastest, it was always the one using (\ref{eq:a-mi}),
but the relative ranking of the other two depended on the model parameters.


\bibliographystyle{agsm}
\bibliography{jors}


\end{document}